\documentclass[12pt]{article}

\usepackage[a4paper,margin=2.5cm]{geometry}

\usepackage{mathrsfs}
\usepackage{latexsym}
\usepackage{amsfonts}
\usepackage{amsmath}
\usepackage{amsthm}
\usepackage{amssymb}
\usepackage[english]{babel}
\input{colordvi}
\usepackage{color}
\usepackage{multirow}
\usepackage{lineno}
\usepackage{array}
\usepackage{multirow}
\usepackage[table]{xcolor}
\usepackage{pstricks}
\usepackage[round]{natbib}   

\usepackage{subfigure}
\usepackage{arydshln} 

\usepackage{epsfig}
\usepackage{algorithm}
\usepackage{algorithmic}

\usepackage{arydshln} 

\title{ {\bf A heuristic algorithm for the \\ Bin Packing Problem with Conflicts \\ on Interval Graphs} }

\author{Tiziano Bacci \footnote{Universit\`a di Roma Tor Vergata, Dipartimento di Ingegneria Informatica e Civile, Via del Politecnico 1, 00133 Roma, Italia, \texttt{bacci.tiziano@gmail.com} } \and Sara Nicoloso\footnote{IASI - CNR, Via dei Taurini 19, 00185 Roma, Italia, \texttt{sara.nicoloso@iasi.cnr.it}} }

\begin{document}

\date{\today}

\maketitle

\begin{abstract}
\noindent 
In this paper we deal with the Bin Packing Problem with Conflicts on interval graphs: given an interval graph, a nonnegative integer weight for each vertex, and a nonnegative integer $B$, find a partition of the vertex set of the graph into $k$ subsets such that the sum of the weights of the vertices assigned to same subset is less than or equal to $B$, two vertices connected by an edge do not belong to the same subset, and $k$ is minimum. We design a heuristic algorithm, and propose a new random interval graph generator which builds interval conflict graphs with desired edge density. We test the algorithm on a huge test bed, and compare the results with existing algorithms.
\end{abstract} 
 
\vfil
\vfil

\section{Introduction}
\label{intro}

In this paper we deal with the Bin Packing problem with Conflicts ($BPPC$) on interval graphs.

$BPPC$, first introduced in a scheduling context (\cite{JO1997}), is defined as follows. Given a graph $G=(V,E)$, a nonnegative integer weight $w_i$ for each vertex $i \in V$, and a nonnegative integer $B$, find a partition of $V$ into $k$ subsets $V_1,\dots,V_k$, such that the sum of the weights of the vertices assigned to a same subset is less than or equal to $B$, two vertices connected by an edge do not belong to the same subset, and $k$ is minimum. Such minimum value of $k$ will be denoted $k_{BPPC}$. The graph $G=(V,E)$ is called {\em conflict graph} and two vertices connected by an edge are said to be {\em in conflict}. $BPPC$ is the union of two well known combinatorial optimization problems, the Bin Packing problem ($BP$) and the Vertex Coloring problem ($VC$), which we now introduce.

The Bin Packing problem ($BP$) is defined as follows. Given a set $V$ of items, a nonnegative integer weight $w_i$ for each item $i \in V$, and a nonnegative integer $B$, find a partition of $V$ into $k$ subsets $V_1,\dots,V_k$, such that the sum of the weights of the items assigned to same subset is less than or equal to $B$ and $k$ is minimum. Such minimum value of $k$ will be denoted $k_{BP}$. Notice that both $BPPC$ and $BP$ can also be defined with rational weights in $[0,1]$ and $B=1$.

The Vertex Coloring problem ($VC$) is defined as follows. Given a graph $G=(V,E)$, find a partition of $V$ into $k$ subsets $V_1,\dots,V_k$, such that two vertices connected by an edge do not belong to the same subset and $k$ is minimum. Throughout the paper $k${\em -coloring} denotes a feasible vertex coloring with $k$ colors. The minimum value of $k$ such that $G$ admits a $k$-coloring is called the {\em chromatic number} $\chi(G)$ of the graph $G$. Notice that each $V_i$ for $i=1,\dots,p$ is an independent set, that is a subset of vertices no two of which are connected by an edge. 

Clearly, $k_{BPPC} \ge \max\{k_{BP},\chi(G)\}$.  

Throughout the paper $n=|V|$, the words {\em vertex} and {\em item} will be used interchangeably, and the subsets of a feasible solution to $BPPC$, $BP$, and $VC$ will be called {\em colors} or {\em bins}.

Since $VC$ on arbitrary graphs and $BP$ are both $NP$-hard (\cite{GJ1978}), $BPPC$ is $NP$-hard too.

Observe that when the edge set $E$ of the graph $G$ is empty, $BPPC$ reduces to $BP$. For increasing $|E|$, the effects of $VC$ on $BPPC$ also increase, because of the increasing number of conflicting pairs of items. Let $t$ be the maximum number of heaviest items whose sum of the weights does not exceed $B$. If $\alpha(G) \le t$, then $BPPC$ reduces to $VC$, where $\alpha(G)$ is the cardinality of a maximum independent set of $G$. As a consequence, when $B \ge \sum_{i \in V} w_i$, $BPPC$ reduces to $VC$.

An application of $BPPC$ is discussed in \cite{CMT1979}, where some flammable, explosive, or toxic substances cannot be placed in the same vehicle. 

\vspace{0.3cm}
\noindent
In this paper we focus on $BPPC$ where $G=(V,E)$ is an interval graph. A graph \linebreak $G=(V,E)$ is an interval graph if every vertex $p \in V$ can be put in one-to-one correspondence with an open intervals $I_p=(l_p,r_p )$ of the real line, and two vertices $p,q \in V$ are connected by edge $(p,q) \in E$ if and only if the corresponding intervals intersect, i.e.~$l_p < r_q$ and $l_q < r_p$. The family of intervals $\mathcal I =\{I_h=(l_h,r_h ), h=1,\dots,n\}$ is called an {\em interval model} for $G$. Any interval graph admits an interval model. Notice that when the edge set of the graph is empty, any set of $n$ mutually non-intersecting intervals is an interval model for $G$. In what follows, w.l.o.g.~we can assume that $\min \{l_j, j=1,\dots, n\} = 0$, and define $R= \max \{r_j, j=1. \dots, n\}$. It is worth observing that $VC$ is solvable in linear time on interval graphs, nevertheless $BPPC$ with an interval conflict graph remains $NP$-hard.



\vspace{0,5cm}
\noindent
In this paper we design a new heuristic algorithm for $BPPC$ with interval conflict graphs and test it on thousands of instances. 

As far as we know, no tests on instances of $BPPC$ with arbitrary interval conflict graphs were performed in the literature. In fact, \cite{SV2013} realize that the conflict graphs of the benchmark instances by \cite{MIMT2010} are  interval graphs, and not arbitrary graphs: actually, the conflict graphs of these instances are not arbitrary interval graphs, but special ones, namely threshold graphs (see Section \ref{sec:exactresultsTM}).

Since the existing random interval graph generators output graphs with edge density in a very narrow range and we want to test our algorithm on interval graphs with edge density ranging from 0 to 1, we design a new random interval graph generator which outputs interval graphs with desired edge density. We use this generator to generate 4000 interval conflict graphs, each of which is associated with a set of item weights. By varying $B$ in ten different ways, we obtain 40000 $BPPC$ instances and test our algorithm on all of them. We also implement some heuristics which are an adaptation of three classical algorithms for $BP$, as in \cite{MIMT2010}, run them on the same test bed and compare the results with ours.

The paper is organized as follows. In Section \ref{sec:literature} we review the literature on $BPPC$. In Section \ref{sec:heuristicinterval} a new heuristic algorithm  is proposed, while a new random interval graph generator is found in Section \ref{Sec:GenInt}. Experimental results are discussed in Section \ref{sec:resultsIG}, where we compare our results with other existing approaches. In Section \ref{sec:exactresultsTM}, our heuristic is tested over instances taken from the literature. Section \ref{sec:conclusions} concludes. 

\section{Literature review}
\label{sec:literature}

The $BPPC$ is widely discussed in literature. We start by surveying papers where approaches to $BPPC$ with arbitrary conflict graphs are proposed, then we discuss some issue about literature $BPPC$ instances and some papers devoted to $BPPC$ on interval graphs. 

\cite{KKT2003} present a heuristic algorithm which repeatedly create a new bin for a maximal independent subset with weight close to $B$. \cite{GLS2004} propose a lower bound and six heuristics: one is a direct adaptation of the First-Fit Decreasing algorithm by \cite{J1974}, three are based on graph coloring, and two are based on finding large cliques. They also describe a random graph generation scheme. \cite{BW2005} compare their heuristic algorithm with two of the best approaches by \cite{GLS2004}, with the algorithm by \cite{KKT2003}, and with a standard beam search algorithm: the proposed algorithm outperforms all the others, on average. \cite{MG2009} present a lower bound which outperforms those by \cite{GLS2004} and by \cite{MIMT2009} and is based on iterative runs of the lower bound algorithms by \cite{GLS2004}. A Column Generation approach is proposed by \cite{JMSSV2010} and tested on 280 instances with density between 10\%-40\% (the generation scheme is not specified). An exact algorithm based on a set-covering formulation is discussed by \cite{MIMT2010}. The authors propose a very effective but time consuming lower bounds; the upper bounds are obtained by means of fast and good heuristic algorithms which are an adaptation of the classical First-Fit Decreasing, Best-Fit Decreasing, Worst-Fit Decreasing for $BP$; lower and upper bounds are better than those by \cite{GLS2004}. If no optimal solution is found then a population-based metaheuristic is applied and possibly a  Branch-and-Price algorithm is adopted. \cite{KCT2010} improve some lower bounds by \cite{MIMT2010} by applying reduction procedures. \cite{ELGN2011} propose a Branch-and-Price algorithm which is compared with those by \cite{MIMT2010}: results show that neither one outperforms the other. \cite{MR2011} propose seven heuristics: one is an adaptation of the Minimum Bin Slack heuristic by \cite{GH1999}; the others six repeatedly create a new bin by selecting (by means of classical bin packing methods) a subset of items from a maximal independent set previously generated. The experimental results show that these heuristics outperform those by \cite{GLS2004}. \cite{YLW2014} show the effectiveness of an ant colony optimization approach to determine a feasible coloring solution to which an improved First-Fit Decreasing heuristic bin packing procedure is applied. \cite{GI2016} describe an effective Column Generation approach to solve $BPPC$ and other problems to optimality, providing new classes of valid inequalities. \cite{BP2016} present an exact method based on an arc-flow formulation with side constraints. The method builds very strong integer programming models that can be given in input to any state-of-the-art mixed integer programming solver. The algorithm is applied to many classical combinatorial problems and, in particular, all the instances by \cite{MIMT2010} are efficiently solved to optimality.

We remark that the generator by \cite{GLS2004} has been improperly used to generate arbitrary graphs (\cite{BW2005, BP2016, CFVO2015, CKHT2011, CFM2017, ELGN2011, GI2016, ThJoncour2010, JMSSV2010, JOK2015, ThKhanafer2010, KCHT2012, KCT2010, KCT2012, MG2009, MR2011, ThMuritiba2010, MIMT2010, SV2013, YLW2014}). In Section \ref{sec:exactresultsTM} we show that these are not arbitrary graphs but special interval graphs, namely threshold graphs. $BPPC$ on threshold graphs turns out to be easier than on arbitrary interval graphs and arbitrary graphs (see \cite{BN2017TG}). We remark that \cite{MIMT2010} used the generator by \cite{GLS2004} to build publicly available instances (see http://or.dei.unibo.it/library/bin-packing-problem-conflicts) and many of the authors above used them. 

Few papers are devoted to $BPPC$ with interval conflict graphs. \cite{EL2008} present a $\frac{7}{3}$-approximation algorithm. \cite{SV2013} present a generic effective Branch-and-Price algorithm for the $BPPC$ with arbitrary conflict graphs, using a Dynamic Programming algorithm for pricing when the conflict graph is an interval graph. They test their algorithm on the instances by \cite{MIMT2010} and \cite{ELGN2011}, closing all open instances, and on harder instances with an arbitrary conflict graph and a larger number of items per bin. A $2$-approximation algorithm exists if $G$ is a threshold graph. In fact, \cite{M2008} proposes a 2-approximation algorithm for the Minimum Clique Partitioning Problem on a weighted interval graph: given an interval graph with nonnegative vertex weights, find a partition of the vertices into the minimum number of cliques such that the sum of vertex weights in each clique does not exceed a given bound $B$. Since $G$ is a threshold graph, it is also a co-interval graph, hence $\overline{G}$, the complement of $G$, is an interval graph and $BPPC$ on a threshold graph $G$ is equivalent to the Minimum Clique Partitioning Problem on the weighted interval graph $\overline{G}$.

We recall the following problems related to $BPPC$. Given an instance of $BPPC$, let $D_p = \{i \in V: w_i = W_p\}$ be the set of items  whose weight is equal to $W_p$, and let $d_p = |D_p|$. If we say that the items in $D_p$ are {\em items of type $p$} and that $d_p$ is their demand, then $BPPC$ is usually known as Cutting Stock problem ($CS$) (\cite{DIM2016}) and it is formulated with integer variables and not binary ones like $BPPC$. $CS$ arises in industrial contexts and often the number of different types of items is small w.r.t.~the number $n$ of items, while in $BPPC$ and $BP$ it is not. The special case of BPPC when $w_i = 1$ for $i=1,\dots,n$ is the optimization version of Bounded Independent Sets (also known as Mutual Exclusion Scheduling, see \cite{BC1996}), which is $NP$-complete when $G$ is an interval graph and $B \ge 4$ (\cite{BJ1995}). An application of Mutual Exclusion Scheduling on interval conflict graphs is described in \cite{G2009}. In \cite{GBSCC2008} also item-bin conflicts are considered. When the number of bin is fixed, \cite{KL2015} minimize the weight of the heaviest bin, while \cite{KCHT2012} minimize the number of violated conflicts. 

\section{Heuristic algorithm}
\label{sec:heuristicinterval}

Since $BPPC$ is the union of $BP$ and $VC$, we believe that in order to design effective heuristic algorithms one has to adopt one of the following two approaches, where \linebreak $P \in \{VC,BP\}$: modify an algorithm designed for a problem $P$ to directly obtain a feasible solution for $BPPC$, or determine a feasible solution for problem $P$, first, then modify it to obtain a feasible solution for $BPPC$.

The first approach is used by the algorithms $U_{FF(\alpha)}$, $U_{BF(\alpha)}$, and $U_{WF(\alpha)}$ by \cite{MIMT2010}: an algorithm designed for $BP$ is modified in such a way that an item $p$ is not assigned to a bin containing an item $q$ in conflict with $p$. The first approach is applied also by \cite{GLS2004}: an algorithm designed for $VC$ is modified in such a way that when a new color $S$ is created and its weight exceeds $B$, then a suitable subset $S' \subset S$ is determined such that $\sum_{i \in S \setminus S'} w_i \le B$, and removed from $S$.

Here we adopt the second approach, as $VC$ on interval graphs is solvable in linear time. In particular we construct an optimum feasible vertex coloring solution and modify it with a local search approach to obtain a feasible $BPPC$ solution.

\vspace{0,3 cm}
\noindent   
In the algorithm, we will make use of an interval model of $G$ and of the following notations or definitions. 
\begin{itemize}
\item  $\lambda = \max \{LB_{BP}, \omega(G)\}$, a lower bound for $BPPC$, where $LB_{BP}$ denotes a lower bound for the Bin Packing problem underlying the given $BPPC$ and $\omega(G)$ denotes the size of a maximum clique of $G$ (recall, in fact, that $\omega(G)=\chi(G)$ as $G$ is an interval graph); 
\item $\mathcal C$, the leftmost subset of $\omega(G)$ mutually intersecting intervals, i.e.~the leftmost maximum clique;
\item $\pi = \max \{l_j, I_j\in \mathcal{C} \}$, the leftmost coordinate belonging to $\omega(G)$ intervals; 
\item $\mathcal I_{left}$, the set of intervals whose right endpoint lays on the left of $\pi$;
\item  $\mathcal I_{right}$, the set of intervals whose left endpoint lays on the right of $\pi$;
\item $W^{est}(V_i)$, the {\em estimated weight} of the empty space on the right of the interval belonging to $V_i \cap \mathcal C$, if any, or on the right of $\pi$: assuming that the interval weights are uniformly distributed among the $R \times \lambda$ unit segments, let $\mu=\frac {1}{R\lambda} \sum_{j=1, \dots, n} w_j$ be the average weight of a unit segment, denote by $I_j=(l_j,r_j)$ the (unique) interval in $\mathcal C$ which belongs also to $V_i$ for $i=1,\dots,\omega(G)$, then $R_i= r_j$ for $i=1,\dots,\omega(G)$, and $R_i= \pi$  for $i = \omega(G) + 1,\dots,\lambda$, and $W^{est}(V_i)= \mu (R-R_i)$;
\item $z$ and $V_1$, \dots , $V_z$, the number of subsets and the subsets in the current (possibly infeasible) solution, respectively;
\item $W(X)$, the weight of a subset $X$ of intervals, i.e.~the sum of the weights of the intervals belonging to $X \subseteq V$; if $W(X) > B$ the subset $X$ is {\em heavy}, otherwise it is {\em light};
\item $Tail(V_i,\rho) \subseteq V_i$, where $\rho \in [0,R]$ denotes a coordinate, a subset of intervals of $V_i$ which we define only for those $i$ such that $\nexists I_j \in V_i : \; l_j <\rho< r_j$; when defined, $Tail(V_i,\rho) = \{I_j \in V_i : \; l_j \ge \rho\}$.
\item a subset $V_i$  {\em non-conflicting} w.r.t.~interval $I_j$, i.e. an independent subset $V_i$ such that $V_i \cup \{I_j\}$ is an independent subset too. 
\end{itemize}

\vspace{0,3cm}
\noindent
The algorithm consists of two phases.
 
In the first phase the algorithm constructs a $\lambda$-coloring $\{ V_1,$\dots,$V_\lambda \}$ of $G$ working on the chosen interval model of $G$. In particular, among all the feasible $\lambda$-colorings of $G$, the algorithm finds a coloring where the weight of the lightest color is as large as possible, as we now describe.

It starts by assigning each interval of the leftmost maximum clique $\mathcal C$ to a different subset, then it assigns the intervals on the left of $\pi$, and finally those on the right of $\pi$. On the left of $\pi$, the algorithm repeatedly assigns an unassigned interval with rightmost right endpoint to the color $V_i$ with minimum (current) weight $W(V_i) + W^{est}(V_i)$. On the right of $\pi$, the algorithm repeatedly assigns an unassigned interval with leftmost left endpoint to the color $V_i$ with smallest (current) weight $W(V_i)$.

The algorithm ends the first phase with a feasible $\lambda$-coloring $\{ V_1,$\dots,$V_\lambda \}$ of $G$. If \linebreak $W(V_i) \le B$ for $i=1,\dots,\lambda$, the partition $\{ V_1,$\dots,$V_\lambda \}$ is a feasible solution (of value $\lambda$) for $BPPC$. Since its value equals the lower bound, $\{ V_1,$\dots,$V_\lambda \}$ is also an optimum solution for $BPPC$. If this is not the case, the algorithm proceeds with the second phase.
 
\vspace{0.3cm}
\noindent
In the second phase the algorithm repeatedly selects a heaviest subset $V_g$ and suitably modifies it to get a light subset. This is accomplished in two different ways: the {\sc tail-exchange}, 
and the {\sc insertion}
.
 
The {\sc tail-exchange} w.r.t.~a coordinate $\rho$ between the chosen heavy subset $V_g$ and a light subset $V_h$ with minimum $W(Tail(V_h,\rho))$, consists of exchanging $Tail(V_g,\rho)$ with $Tail(V_h,\rho)$. It can be done iff the following three conditions are verified: both $Tail(V_g,\rho)$ and $Tail(V_h,\rho)$ are defined, $W(Tail(V_h,\rho)) < W(Tail(V_g,\rho))$, and the resulting $V_h$ keeps being light. Notice that the weight of $V_g$ after the exchange is decreased. Precisely, the algorithm finds the leftmost coordinate $\rho \ge \min \{r_s: I_s \in V_g\}$ such that there exists a subset $V_h$ allowing a {\sc tail-exchange} operation with $V_g$. If the resulting $V_g$ is still heavy, the algorithm finds the next $\rho$ with the same properties and repeats this step again, stopping as soon as $\rho \ge R$ or the current $V_g$ is light. If $V_g$ is still heavy the algorithm tries to apply the {\sc insertion}.

In an {\sc insertion}, an interval $I_j \in V_g$ which minimizes $|W(V_g) - w_j - B|$ is selected, and inserted into a light non-conflicting $V_h$ such that the resulting $V_h$ is light and $W(V_h) + w_j$ is maximum, if any. Otherwise $I_j$ is inserted into a heavy non-conflicting $V_h$ with minimum weight, if any. The first time that no {\sc insertion} is performed w.r.t.~$V_g$, then a new subset is created and $I_j$ is inserted into it. 

\vspace{0.3cm}
\noindent
The algorithm description follows. By {\sc tail-exchange step} ({\sc insertion step}, respectively) we mean the repeated application of a {\sc tail-exchange}  ({\sc insertion}, respectively). We recall that, when $\omega(G) =1$ (i.e.~the edge density is zero), any set of $n$ mutually non-intersecting intervals is an interval model of $G$, and $BPPC$ reduces to $BP$.

\vspace{1cm}
{\small
\noindent
ALGORITHM BN
\vspace{0.3cm}
{\tt

\noindent
{\em Input}:~an interval model for the graph $G$, $w_i \in \mathbb Z_+$ $\forall i \in V$, $B \in \mathbb Z_+$.

\noindent
{\em Output}:~$z$ and a feasible partition $\{V_1,\dots,V_z\}$ of $V$.
\vspace{0.3cm}
\begin{itemize}
\setlength{\parskip}{-3pt}
\item {\sc Phase I}
\begin{itemize}
\setlength{\parskip}{-3pt}
\item [] Define $z:= \lambda$;
\item [] Define $z$ empty sets $V_1$, \dots , $V_z$;
\item [] Assign each interval of $\mathcal C$ to a different subset;
\item [] Let $I_j \in \mathcal{I}_{left}$ be an interval with righmost right endpoint,
\begin{itemize}
\setlength{\parskip}{-3pt}
\item[] assign $I_j$ to a non-conflicting subset $V_i$ with minimum $W(V_i) + W^{est}(V_i)$,
\item[] remove $I_j$ from $\mathcal{I}_{left}$, and repeat until $\mathcal{I}_{left}$ is empty;
\end{itemize}
\item [] Let $I_j \in \mathcal{I}_{right}$ be an interval with leftmost left endpoint,
\begin{itemize}
\setlength{\parskip}{-3pt}
\item[] assign $I_j$ to a non-conflicting subset $V_i$ with minimum $W(V_i)$,
\item[] remove $I_j$ from $\mathcal{I}_{right}$, and repeat until $\mathcal{I}_{right}$ is empty;
\end{itemize}
\end{itemize}
\item {\sc Phase II}
\begin{itemize}
\setlength{\parskip}{-3pt}
\item [] While $V_1, \dots , V_z$ is infeasible do
\begin{itemize} 
\setlength{\parskip}{-3pt}
\item [] Let $V_g$ be a subset with maximum weight; 
\item []  $\rho := \min \{r_s: I_s \in V_g\}$;
\item [] {\sc tail-exchange step}: 
\begin{itemize}
\setlength{\parskip}{-3pt}
\item[] While $W(V_g) > B$ and  $\rho < R$ do 
\item[] let $V_h$, $h \neq g$, with minimum $W(Tail(V_h,\rho))$ be a light subset  
\item[] such that $W(Tail(V_h,\rho)) < W(Tail(V_g,\rho))$ and
\item[] $W(V_h) - W(Tail(V_h,\rho)) + W(Tail(V_g,\rho)) \le B$, if any;
\item[] set $V_g := V_g \setminus Tail(V_g,\rho) \cup Tail(V_h,\rho)$ 
\item[] and $V_h := V_h \setminus Tail(V_h,\rho) \cup Tail(V_g,\rho)$
\item[] $\rho = \rho +1$;
\end{itemize}
\item [] new\_subset:={\sc false};  
\item [] {\sc insertion step}:
\begin{itemize}
\setlength{\parskip}{-3pt}
\item [] While $W(V_g) > B$ do 
\item[] let $I_j \in V_g$ be an interval with minimum $|W(V_g) - B - w_j|$
\item[] If there exists a $V_h$, $h \neq g$, non-conflicting w.r.t.~$I_j$, 
\item[] such that $W(V_h) + w_j$ is maximum and $\le B$,  
\item[] then remove $I_j$ from $V_g$, and insert it in $V_h$,
\item[] otherwise if there exists a $V_h$, $h \neq g$, non-conflicting 
\item[] w.r.t.~$I_j$, such that $W(V_h)$ is mimimum and $\ge B$, 
\item[] then remove $I_j$ from $V_g$, and insert it in $V_h$, 
\item[] otherwise if new\_subset={\sc false} set $z:=z+1$ and new\_subset:={\sc true}
\item [] then remove $I_j$ from $V_g$ and insert it in $V_z$.
\end{itemize}    
\end{itemize}
\end{itemize}
\end{itemize}
}
}

\vspace{0,3 cm}
\noindent
If the partition at the end of Phase I is feasible then the algorithm terminates, otherwise the algorithm enters Phase II and chooses a heaviest subset $V_g$. If $V_g$ after the {\sc tail-exchange step} is light then the algorithm terminates an iteration of the main {\sc while}-instruction of Phase II and the number of heavy subsets is decreased by one, otherwise the algorithm performs the {\sc insertion step}: if no new subset is created then, again, the algorithm terminates an iteration of the main {\sc while}-instruction of Phase II and the number of the heavy subsets is decreased by one; on the contrary, if a new subset $V_{z+1}$ is created, the algorithm terminates this iteration of the main {\sc while}-instruction either with $V_g$ and $V_{z+1}$ light or with $V_g$ light and $V_{z+1}$ heavy but verifying $W(V_{z+1}) < W(V_g)$. In the former case the number of heavy subsets is decreased by one, in the latter the number of heavy subsets keeps constant but the overall infeasibility is decreased. This discussion show that the algorithm terminates.   

The computational complexity of Phase I is $O(n\log n)$. As for Phase II, the complexity of the {\sc tail-exchange step} is $O(nR)$ and the complexity of the {\sc insertion step} is $O(n^2)$. Since the main {\sc while}-instruction is repeated at most $n$ times and $R \le 2n$ (in fact, w.l.o.g., one can delete all the coordinates which belong to at most one interval), the overall computational complexity of Phase II is $O(n^3)$, and so is the complexity of the entire algorithm.

The algorithm has been tested over thousands of randomly generated instances. Computational results are presented in Section \ref{sec:resultsIG}.

\section{A random interval graph generator}
\label{Sec:GenInt}

\noindent
The easiest way to generate a random interval graph is to randomly choose the endpoints of each of the $n$ intervals and then construct the intersection graph of them (\cite{ThVasileios2005,JSW1990}).  

Another (equivalent) generator is the following. Given a non-negative integer $n$, let $\pi = (\pi_1,\dots,\pi_{2n})$ be a random permutation of $(1,1,2,2,\dots,n,n)$; then for $j=1,\dots,n$ define interval $I_j=(l_j,r_j)$, where $l_j = \min \{k: \pi_k = j, k=1,\dots,2n\}$ and $r_j = \max \{k: \pi_k = j, k=1,\dots,2n\}$ (notice that in $\pi$ there exists exactly two elements of value $j$). 

We generated thousands of set of intervals in both ways. The experimental analysis we conducted shows that the edge density $2|E|/ (n (n-1))$ of almost all the corresponding interval graphs is 60\%-70\%. The $BPPC$ instances in the literature (see, for example, \cite{GLS2004,MIMT2010,SV2013}) are classified by their edge density ranging from $0\%$ to $90\%$. Since, to our knowledge, no random interval graph generator exists which allows to obtain an interval graph with a prescribed edge density, in the present section we define a new random interval graph generator with this property. In fact, the edge density of the intersection graph of a set of intervals in $[0;D]$ depends on the average interval length.  

Our generator accepts in input the number $n$ of intervals and the desired edge density $\delta$ of the corresponding intersection graph, and suitably computes three values, $D$, $\Lambda_{\min}$,  $\Lambda_{\max}$ to obtain an interval graph with desired edge density $\delta$. The output consists of a set of $n$ intervals whose endpoints (integer, w.l.o.g.)~are uniformly distributed in $\{0,\dots,D\}$, and an arbitrary interval $I_j=(l_j,r_j)$ in this set will have  length $\Lambda_j = r_j - l_j$ verifying $\Lambda_{\min} \le \Lambda_j \le \Lambda_{\max}$. Experimental analysis conducted on thousands of graphs shows that our generator generates arbitrary interval graph with expected edge density $\delta$ (the standard deviation of $\delta$ increases for increasing $\delta$ and ranges in [0,0.02]). Let's go into details.

In order to ensure that the $2n$ endpoints of the $n$ intervals have space enough to give a graph with edge density equal to zero (for suitable interval lengths), one has to fix $D \ge 2n$. We tried many different values for $D$ but we did not appreciate any differences, so we decided to set $D = 2.5n$.

The edge density $\delta$ of the resulting interval graph depends on the average interval length $\overline{\Lambda}$: for example when  $\overline{\Lambda} = 1$ then $\delta \simeq 0$, and when  $\overline{\Lambda} > D/2$ then $\delta \simeq 1$. We can determine the equation which ties $\delta$ and $\overline{\Lambda}$. The coordinate $0$ can be chosen as a left endpoint, only, the coordinate $D$ can be chosen as a right endpoint, only, and all the other coordinates can be chosen both as left and right endpoints. Hence the average number of  left endpoints per coordinate is $\frac{n}{D}$, as well as the average number of right endpoints per coordinate, and an interval of length $\Lambda$ intersects $\frac{2n}{D} \Lambda$ intervals, on average. Thus the average degree of a vertex is $\frac{2n}{D} \overline{\Lambda}$, resulting in $\frac{2n} {D} \frac{n} {2} \overline{\Lambda}$ edges of the interval graph. Dividing this quantity by $\frac{n(n-1)}{2}$, one gets that the edge density is $\delta = \frac{2n}{(n-1)D} \overline{\Lambda}$ from which we derive that the average inteval length $\overline{\Lambda}$ has to be set to $\delta D \frac{n - 1}{2n}$ in order to obtain an interval graph with expected edge density $\delta$. 

Recalling that $\overline{\Lambda} = 1/n \sum_{j=1}^n \Lambda_j$ and that $\Lambda_{\min} \le \Lambda_j \le \Lambda_{\max}$ for all $j$, the set of intervals has average length $\overline{\Lambda}$ if one suitably chooses $\Lambda_{\min}$ and $\Lambda_{\max}$. We decided to randomly choose $\Lambda_{\min}$ in a suitable range which we discuss in a while; $\Lambda_{\max}$ is consequently determined. We reason as follows.

\begin{figure}[htbp] 
  \begin{center}
    \begin{pspicture}(1,0)(15,5)

\psline{->}(1.75,-0.25)(1.75,4.3)
\rput (1.2,4) {{\scriptsize {\sc eil}$(l_j)$}}

\psline{->}(1.5,0)(15.5,0)
\rput (15.3,-0.25) {{\scriptsize $l_j$}}

\rput (1.1,3.2) {{\scriptsize $\Lambda_{\max}$}}
\psline(1.65,3.2)(1.85,3.2)
\rput (1.1,1) {{\scriptsize $\Lambda_{\min}$}}
\psline(1.65,1)(1.85,1)
\rput (0.8,2.1) {{\scriptsize $\frac{\Lambda_{\max} + \Lambda_{\min}}{2}$}}
\psline(1.65,2.1)(1.85,2.1)

\rput (2.95,-0.3) {{\scriptsize $\Lambda_{\min}$}}
\psline(2.95,-0.1)(2.95,0.1)
\rput (5.2,-0.3) {{\scriptsize $\Lambda_{\max}$}}
\psline(5.2,-0.1)(5.2,0.1)

\rput (13.5,-0.3) {{\scriptsize $D$$-$$\Lambda_{\min}$}}
\psline(13.5,-0.1)(13.5,0.1)
\rput (11.3,-0.3) {{\scriptsize $D$$-$$\Lambda_{\max}$}}
\psline(11.3,-0.1)(11.3,0.1)

\rput (14.5,-0.3) {{\scriptsize $D$}}
\psline(14.5,-0.1)(14.5,0.1)

\psline(1.75,2.1)(11.3,2.1)

\psline(11.3,2.1)(13.5,1)

\psline[linestyle=dashed,dash=2pt 2pt](13.5,0)(13.5,1)
\psline[linestyle=dashed,dash=2pt 2pt](11.3,2.1)(11.3,0)
\psline[linestyle=dashed,dash=2pt 2pt](1.75,1)(13.5,1)

   \end{pspicture}
    \end{center}
    \vspace{0.5cm}
  \caption{The expected interval length {\sc eil}$(l_j)$ when $l_j$ is randomly chosen}
  \label{fig:lambdamedio}
\end{figure}
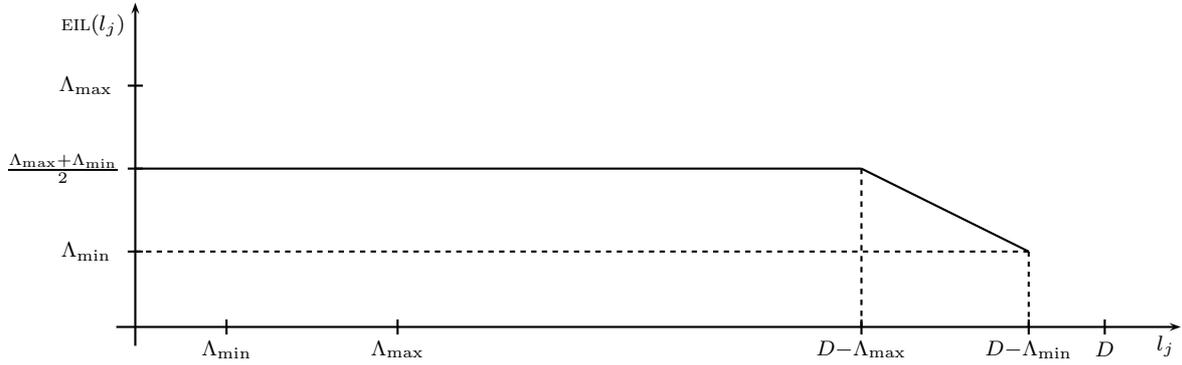

By flipping a coin, the generator randomly chooses the left endpoint $l_j$ of interval $I_j$ first, or the right endpoint $r_j$. Assume that the left endpoint $l_j$ is randomly chosen first. In order to ensure that the length $\Lambda_j$ verifies $\Lambda_j \ge \Lambda_{\min}$, $l_j$ has to be chosen in $\{0,\dots,D - \Lambda_{\min}\}$. Refer to Figure \ref{fig:lambdamedio}. If $l_j \in \{0,\dots,D - \Lambda_{\max}\}$ then $r_j$ can be randomly chosen in $\{l_j + \Lambda_{\min},\dots,l_j + \Lambda_{\max}\}$. In this case the expected interval length {\sc eil}$(l_j)$ is $(\Lambda_{\min} + \Lambda_{\max})/2$. If $l_j \in \{D - \Lambda_{\max} + 1,\dots,D - \Lambda_{\min}\}$ then $r_j$ can be randomly chosen only in $\{l_j + \Lambda_{\min},\dots,D\}$ and the {\sc eil}$(l_j)$ is a linear function of $l_j$, namely $(D + \Lambda_{\min} - l_j)/2$. Thus the average interval length  $\overline{\Lambda}'$ when $l_j$ is chosen first can be computed dividing the entire area underlying the drawn function by its width $D - \Lambda_{\min}$ as follows:
$$\overline{\Lambda}' =  \frac{1}{D - \Lambda_{\min}}\big [\big (D - \Lambda_{\max} \big) \frac{\Lambda_{\min} + \Lambda_{\max}}{2} +  \big (\frac{\Lambda_{\max} - \Lambda_{\min}}{2} \big ) \big (\frac{\Lambda_{\min} + \Lambda_{\max}}{2} + \Lambda_{\min} \big ) \big ] =$$ $$= \frac{1}{4(D - \Lambda_{\min})} [-\Lambda_{\max}^2 - 3\Lambda_{\min}^2 + 2D(\Lambda_{\max}+\Lambda_{\min})] $$
\noindent
Assume now that the right endpoint $r_j$ of interval $I_j$ is randomly chosen first. By similar arguments, $r_j$ has to be chosen in $\{\Lambda_{\min},\dots,D\}$ and $l_j \in \{\max\{r_j - \Lambda_{\max};0\},\dots,r_j-\Lambda_{\min} \}$ and the average interval length $\overline{\Lambda}''$ when $r_j$ is chosen first is equal to $\overline{\Lambda}'$ (the graph of {\sc eil}$(r_j)$ can be obtained by horizontally flipping the graph of Figure \ref{fig:lambdamedio} along the vertical axis $D/2$). 

Since, on average, one half of the intervals are generated by randomly choosing the left endpoint first, and one half by randomly choosing the right endpoint first, the overall average interval length is

$$\overline{\Lambda}= \frac{\overline{\Lambda}' + \overline{\Lambda}''}{2} = \frac{1}{4(D - \Lambda_{\min})}  [-\Lambda_{\max}^2 - 3\Lambda_{\min}^2 + 2D(\Lambda_{\max}+\Lambda_{\min})] $$

\noindent
Given $\Lambda_{\min}$ and $\overline{\Lambda}$, we can use this formula to determine $\Lambda_{\max}$. Since $\Lambda_{\max} >  D$ makes no sense, we get
$$\Lambda_{\max} =  D - \sqrt{D^2 - 4 D \overline{\Lambda} + 2 D \Lambda_{\min} + 4 \overline{\Lambda} \Lambda_{\min} - 3 \Lambda_{\min}^2}$$
\noindent
The non-negativity of the argument of the square root requires $\frac{4 \; \overline{\Lambda} - D}{3} \le \Lambda_{\min} \le D$. On the other hand, clearly, $1 \le \Lambda_{\min} \le \overline{\Lambda}$. Hence 
$\max\{1; \big \lceil \frac{4 \; \overline{\Lambda} - D}{3} \big \rceil \} \leq \Lambda_{\min} \leq \overline{\Lambda}$, as $ \overline{\Lambda} \le D$.

The proposed generator fixes $D$ and computes $\overline{\Lambda}$ as discussed. Then for each interval $I_j$, it randomly chooses $\Lambda_{\min} \in \left \{\max\{1; \big \lceil \frac{4 \; \overline{\Lambda} - D}{3} \big \rceil \} \; ,\dots, \; \overline{\Lambda} \right\}$, computes the corresponding $\Lambda_{\max}$, randomly chooses whether to generate $r_j$ first or $l_j$ and, in both cases, randomly generates a suitable $\Lambda_j$, and consequently fixes the other endpoint.  

\vspace{1 cm}
{\small
\noindent
RANDOM INTERVAL GRAPH GENERATOR
\vspace{0.3cm}
{\tt

\noindent
{\em Input}:~$n \in \mathbb Z_+$ and $\delta \in [0,1]$

\noindent
{\em Output}:~a set of $n$ intervals whose intersection graph has expected edge density $\delta$.\vspace{0.3cm}
\begin{itemize}
\setlength{\parskip}{-3pt}
\item [] $D := 2.5 \; n$;
\item [] $\overline{\Lambda} := \delta D \frac{n - 1}{2n}$;
\item [] For $j=1,\dots,n$
 \begin{itemize}
\setlength{\parskip}{-3pt}
	\item [] randomly choose $\Lambda_{\min} \in \{ \; \max\{1; \big \lceil \frac{4 \; \overline{\Lambda} - D}{3} \big \rceil \} \; ,\dots, \; \overline{\Lambda} \; \}$
	\item [] $\Lambda_{\max} :=  D - \sqrt{D^2 - 4 D \overline{\Lambda} + 2 D \Lambda_{\min} + 4 \overline{\Lambda} \Lambda_{\min} - 3 \Lambda_{\min}^2}$;
 \item [] flip a $coin$;
 \item [] if $coin$ = {\sc head} then randomly choose $r_j \in \{\Lambda_{\min},\dots,D\}$ and 
 \begin{itemize}
\setlength{\parskip}{-3pt}
\item[] $\Lambda_j \in \big \{\Lambda_{\min},\dots,\min \{ r_j ; \Lambda_{\max}\} \big \}$, and set $l_j := r_j - \Lambda_j$;
 \end{itemize}
 \item [] if $coin$ = {\sc tail} then randomly choose $l_j \in \{0,\dots,D-\Lambda_{\min}\}$ and 
 \begin{itemize}
\setlength{\parskip}{-3pt}
\item [] $\Lambda_j \in \big \{ \Lambda_{\min},\dots,\min \{ D - l_j ; \Lambda_{\max} \} \big \}$, and set $r_j := l_j + \Lambda_j$;
\end{itemize}
\end{itemize}
\end{itemize}
}
}

\section{Computational results}
\label{sec:resultsIG}

\noindent
In the present section we discuss the results obtained by solving thousands of instances.

The test bed was generated as we now describe. By $TI(n,B,\Delta)$ we denote a set of 100 randomly generated instances of $BPPC$ with $n$ items, weights uniformly distributed in $[20,100]$ (as in \cite{F96}), bound $B$, and interval conflict graph with expected edge density $\Delta$. When $\Delta > 0$ we repeatedly run the random interval graph generator described in Section \ref{Sec:GenInt} and we selected 100 sets of $n$ intervals whose intersection graph had edge density $\delta \in [\Delta - 0.02; \Delta + 0.02]$. When $\Delta = 0$  we defined the set $\mathcal I =\{I_h=(h,h+1 ), h=0,\dots,n-1\}$ of $n$ mutually non-intersecting intervals (in this case $BPPC$ reduces to $BP$). In particular, we chose $n\in\{120,250,500,1000\}$, $B \in \{120,$ $150,$ $180,$ $210,$ $240,$ $270,$ $300,$ $330,$ $360,$ $390\}$, and $\Delta \in \{0,$ $0.1,$ $0.2,$ $0.3,$ $0.4,$ $0.5,$ $0.6,$ $0.7,$ $0.8,$ $0.9\}$. Totally we built 40000 instances. 

Notice that on these instances the number of different weights is $100-20+1=81$. Hence, every weight is expected to appear in $n/81$ copies. For increasing $n$ the underlying (classical) Bin Packing recalls a Cutting Stock (see Section \ref{intro}).

We compare the computational results obtained by applying the heuristic algorithm $BN$ proposed in Section \ref{sec:heuristicinterval} and an adaptation to $BPPC$ of the classical heuristic algorithms First-Fit Decreasing, Best-Fit Decreasing, Worst-Fit Decreasing for $BP$ (\cite{J1974}), as described in \cite{MIMT2010}. In particular, these adaptations, $U_{FF(\alpha)}$, $U_{BF(\alpha)}$, and $U_{WF(\alpha)}$ (we shall call them algorithms $M$), consider an extended conflict graph $G_w$, obtained by adding to $G$ an edge for each pair of vertices $i,j$ with $w_i + w_j > B$, and consider  vertex weights $w_i^s$ defined as follows: $w_i^s = \alpha (w_i / \overline{w}) +$ $(1 - \alpha) (\deg(i) / \overline{\deg})$, for $i = 1,2,\dots,n$, where $\alpha \in \{0,0.1,\dots,1\}$, $\deg(i)$ is the degree of vertex $i$ in $G_w$, and $\overline{w}$ and $\overline{\deg}$ are the average weight of the vertices and their average degree in $G_w$, respectively. Notice that $G_w \setminus G$ is a threshold graph (see Section \ref{sec:exactresultsTM}) and its edge density $\delta'$ decreases for increasing $B$. For example, when the weights are uniformly distributed in $[20,100]$ as in our test bed, $\delta' = 0.5$ when $B = 120$, $\delta' = 0.18$ when $B = 150$, $\delta' = 0.03$ when $B = 180$, and $\delta' = 0$ when $B \ge 200$.

All the algorithms were coded in C{\small ++} and run on an Intel Core i7-3632QM 2.20GHz $\times$ 8 (up to 3.2 GHz with turbo boost) with 16 GB RAM under a Linux operating system. 

Let $S$ be an instance of $BPPC$ with an interval conflict graph $G$. Given \linebreak $\alpha \in \{0,0.1,\dots,1\}$, let $u_{x(\alpha)}(S)$ be the value of the solution output by algorithm $U_{x(\alpha)}$ on $S$, for $x \in \{FF, BF, WF\}$. By $u^M(S)= \min \{u_{x(\alpha)}(S), \alpha \in \{0,0.1,\dots,1\}, x \in \{FF, BF, WF\}\}$ we denote the minimum among all the 33 values of the (feasible) solutions output by algorithms $M$ on $S$. By $u^{BN}(S)$ we denote the value of the (feasible) solution output by algorithm $BN$ on $S$. 

We tested each algorithm on all the 40000 instances of our test bed. To evaluate the performances of the algorithms we define $LB_{BPPC}(S) = \max \left \{ \big \lceil \sum_{i \in V}  w_i / B \big \rceil ; \chi(G) \right \}$, a lower bound on the value of an  optimum solution of  $BPPC$ on instance $S$. The detailed results obtained for $n=1000$ are presented in Table  \ref{BNvsM1000TI}; the results for $n \in \{120,$ $250,$ $500,$ $1000\}$ are summarized in Table \ref{BNvsMTot120to1000TI}.

In Table \ref{BNvsM1000TI}, where rows are indexed by $\Delta$ and columns by $B$, we compare the results obtained by applying the algorithms on the 10000 instances $TI(1000,B,\Delta)$. In each cell there are six values, each one averaged over the corresponding 100 instances: M=LB (BN=LB, respectively)~is the percentage of instances $S$ where $u^M(S) = LB_{BPPC}(S)$ ($u^{BN}(S) = LB_{BPPC}(S)$, respectively), i.e. the percentage of instances where $LB_{BPPC}(S)$ allows to certify that the corresponding algorithm found an optimum solution; M$<$BN (BN$<$M, respectively)~is the percentage of instances $S$ where \linebreak $u^M(S) < u^{BN}(S)$ ($u^{BN}(S) < u^{M}(S)$, respectively) (notice that the complement to 100\% of the sum of the last two values is the percentage of instances where $u^M(S)=u^{BN}(S)$); \; Gap\_M (Gap\_BN, respectively)~is the gap $\frac{u^M(S) - LB_{BPPC}(S)}{LB_{BPPC}(S)}$ ($\frac{u^{BN}(S) - LB_{BPPC}(S)}{LB_{BPPC}(S)}$, respectively). A light grey indicates the algorithm which outperforms the other one w.r.t.~the corresponding data. If in a cell the value X=LB is 100\%, then all the data of algorithm X are colored in light cyan, for X $\in \{$M,BN$\}$.

The light cyan cells in Table \ref{BNvsM1000TI} show that algorithm $BN$ solves to optimality all the 100 instances of 36 out of 100 cells, while algorithms $M$ solve to optimality all the 100 instances of 7 out of 100 cells. 

Surprisingly algorithm $BN$ shows a very nice behaviour for $\Delta=0$ and $B \ge 180$, that is to say on instances of Bin Packing (without conflicts). Algorithms $M$ show better results when $B \in \{120,150\}$: we suspect that it is because these algorithms use the informations of the extended graph. As a general result, algorithm $BN$ has a better behaviour than $M$ when $B \ge 180$, except when $\Delta = 0.9$ where the two algorithms are essentially equivalent.

\begin{table}

\begin{center}

{\scriptsize
\begin{tabular}{| c | c | r | c | c | c | c | c | c | c | c | c | c | }

\cline{4-13}
\multicolumn{1}{c}{} & \multicolumn{1}{c}{}                     &                         & \multicolumn{10}{c|}{$B$} \\                                                                                                                                                                                                                                                                                                                                                                                                              

\cline{4-13}
 \multicolumn{1}{c}{} &	\multicolumn{1}{c}{}  &  & 120 & 150 & 180 & 210 & 240 & 270 & 300 & 330 & 360 & 390 \\ \hline
 \multicolumn{1}{|c|}{\multirow{60}{*}{$\Delta$}} &	\multicolumn{1}{c|}{\multirow{6}{*}{0}} & M=LB & 0\% & 0\% & 0\% & 0\% & 0\% & 0\% & 0\% & 0\% & 0\% & 0\% \\
&	 & M$<$BN & \cellcolor[HTML]{D6D5D5} 100\% & \cellcolor[HTML]{D6D5D5} 45\% & 0\% & 0\% & 0\% & 0\% & 0\% & 0\% & 0\% & 0\% \\
&	 	 & Gap\_M & \cellcolor[HTML]{D6D5D5} 1.68\% & \cellcolor[HTML]{D6D5D5} 1.19\% & 0.89\% & 3.27\% & 2.94\% & 2.59\% & 2.34\% & 2.13\% & 1.98\% & 1.78\% \\
&	 & BN=LB & 0\% & 0\% & 0\% & \cellcolor[HTML]{D6D5D5} 19\% & \cellcolor[HTML]{D6D5D5} 54\% & \cellcolor[HTML]{D6D5D5} 78\% & \cellcolor[HTML]{D6D5D5} 76\% & \cellcolor[HTML]{D6D5D5} 81\% & \cellcolor[HTML]{D6D5D5} 88\% & \cellcolor[HTML]{D6D5D5} 87\% \\
&	 & BN$<$M & 0\% & 24\% & \cellcolor[HTML]{D6D5D5} 97\% & \cellcolor[HTML]{D6D5D5} 100\% & \cellcolor[HTML]{D6D5D5} 100\% & \cellcolor[HTML]{D6D5D5} 100\% & \cellcolor[HTML]{D6D5D5} 100\% & \cellcolor[HTML]{D6D5D5} 100\% & \cellcolor[HTML]{D6D5D5} 100\% & \cellcolor[HTML]{D6D5D5} 100\% \\
&	 & Gap\_BN & 7.76\% & 1.4\% & \cellcolor[HTML]{D6D5D5} 0.42\% & \cellcolor[HTML]{D6D5D5} 0.29\% & \cellcolor[HTML]{D6D5D5} 0.18\% & \cellcolor[HTML]{D6D5D5} 0.1\% & \cellcolor[HTML]{D6D5D5} 0.12\% & \cellcolor[HTML]{D6D5D5} 0.1\% & \cellcolor[HTML]{D6D5D5} 0.07\% & \cellcolor[HTML]{D6D5D5} 0.08\% \\ \cline{2-13}

&	\multicolumn{1}{c|}{\multirow{6}{*}{0.1}} & M=LB & 0\% & 0\% & 0\% & 0\% & 0\% & 0\% & 0\% & 0\% & 0\% & 0\% \\ 
&	 & M$<$BN & \cellcolor[HTML]{D6D5D5} 100\% & \cellcolor[HTML]{D6D5D5} 51\% & 1\% & 0\% & 0\% & 0\% & 0\% & 0\% & 0\% & 0\% \\ 
&	 & Gap\_M & \cellcolor[HTML]{D6D5D5} 1.68\% & \cellcolor[HTML]{D6D5D5} 1.26\% & 1.04\% & 1.61\% & 1.2\% & 1.38\% & 1.6\% & 1.67\% & 1.91\% & 2.31\% \\ 
&	 & BN=LB & 0\% & 0\% & 0\% & \cellcolor[HTML]{D6D5D5} 7\% & \cellcolor[HTML]{D6D5D5} 32\% & \cellcolor[HTML]{D6D5D5} 59\% & \cellcolor[HTML]{D6D5D5} 68\% & \cellcolor[HTML]{D6D5D5} 72\% & \cellcolor[HTML]{D6D5D5} 78\% & \cellcolor[HTML]{D6D5D5} 79\% \\ 
&	 & BN$<$M & 0\% & 25\% & \cellcolor[HTML]{D6D5D5} 90\% & \cellcolor[HTML]{D6D5D5} 100\% & \cellcolor[HTML]{D6D5D5} 100\% & \cellcolor[HTML]{D6D5D5} 100\% & \cellcolor[HTML]{D6D5D5} 100\% & \cellcolor[HTML]{D6D5D5} 100\% & \cellcolor[HTML]{D6D5D5} 100\% & \cellcolor[HTML]{D6D5D5} 100\% \\ 
&	 & Gap\_BN & 7.86\% & 1.5\% & \cellcolor[HTML]{D6D5D5} 0.57\% & \cellcolor[HTML]{D6D5D5} 0.39\% & \cellcolor[HTML]{D6D5D5} 0.28\% & \cellcolor[HTML]{D6D5D5} 0.18\% & \cellcolor[HTML]{D6D5D5} 0.16\% & \cellcolor[HTML]{D6D5D5} 0.15\% & \cellcolor[HTML]{D6D5D5} 0.13\% & \cellcolor[HTML]{D6D5D5} 0.14\% \\ \cline{2-13}
&	\multicolumn{1}{c|}{\multirow{6}{*}{0.2}} & M=LB & 0\% & 0\% & 0\% & 0\% & 0\% & 0\% & 0\% & 0\% & 0\% & 0\% \\ 
&	 & M$<$BN & \cellcolor[HTML]{D6D5D5} 100\% & \cellcolor[HTML]{D6D5D5} 45\% & 0\% & 0\% & 0\% & 0\% & 0\% & 0\% & 0\% & 0\% \\ 
&	 & Gap\_M & \cellcolor[HTML]{D6D5D5} 1.69\% & \cellcolor[HTML]{D6D5D5} 1.61\% & 1.75\% & 2.51\% & 3.25\% & 4.1\% & 4.92\% & 6.32\% & 7.65\% & 8.4\% \\ 
&	 & BN=LB & 0\% & 0\% & 0\% & \cellcolor[HTML]{D6D5D5} 1\% & \cellcolor[HTML]{D6D5D5} 7\% & \cellcolor[HTML]{D6D5D5} 26\% & \cellcolor[HTML]{D6D5D5} 23\% & \cellcolor[HTML]{D6D5D5} 19\% & \cellcolor[HTML]{D6D5D5} 15\% & \cellcolor[HTML]{D6D5D5} 10\% \\ 
&	 & BN$<$M & 0\% & 25\% & \cellcolor[HTML]{D6D5D5} 98\% & \cellcolor[HTML]{D6D5D5} 100\% & \cellcolor[HTML]{D6D5D5} 100\% & \cellcolor[HTML]{D6D5D5} 100\% & \cellcolor[HTML]{D6D5D5} 100\% & \cellcolor[HTML]{D6D5D5} 100\% & \cellcolor[HTML]{D6D5D5} 100\% & \cellcolor[HTML]{D6D5D5} 100\% \\ 
&	 & Gap\_BN & 7.91\% & 1.85\% & \cellcolor[HTML]{D6D5D5} 0.72\% & \cellcolor[HTML]{D6D5D5} 0.54\% & \cellcolor[HTML]{D6D5D5} 0.42\% & \cellcolor[HTML]{D6D5D5} 0.35\% & \cellcolor[HTML]{D6D5D5} 0.39\% & \cellcolor[HTML]{D6D5D5} 0.47\% & \cellcolor[HTML]{D6D5D5} 0.56\% & \cellcolor[HTML]{D6D5D5} 0.67\% \\ \cline{2-13}
&	\multicolumn{1}{c|}{\multirow{6}{*}{0.3}} & M=LB & 0\% & 0\% & 0\% & 0\% & 0\% & 0\% & 0\% & 0\% & 0\% & 0\% \\ 
&	 & M$<$BN & \cellcolor[HTML]{D6D5D5} 100\% & 33\% & 0\% & 0\% & 0\% & 0\% & 0\% & 0\% & 0\% & 0\% \\ 
&	 & Gap\_M & \cellcolor[HTML]{D6D5D5} 1.74\% & 2.61\% & 3.63\% & 4.83\% & 7.49\% & 8.29\% & 7.18\% & 5.71\% & 5.31\% & 5.08\% \\ 
&	 & BN=LB & 0\% & 0\% & 0\% & 0\% & 0\% & 0\% & \cellcolor[HTML]{D6D5D5} 7\% & \cellcolor[HTML]{D6D5D5} 57\% & \cellcolor[HTML]{D6D5D5} 97\% & \cellcolor[HTML]{E0FFFF}  100\% \\ 
&	 & BN$<$M & 0\% & \cellcolor[HTML]{D6D5D5} 61\% & \cellcolor[HTML]{D6D5D5} 100\% & \cellcolor[HTML]{D6D5D5} 99\% & \cellcolor[HTML]{D6D5D5} 100\% & \cellcolor[HTML]{D6D5D5} 100\% & \cellcolor[HTML]{D6D5D5} 100\% & \cellcolor[HTML]{D6D5D5} 100\% & \cellcolor[HTML]{D6D5D5} 100\% & \cellcolor[HTML]{E0FFFF} 100\% \\ 
&	 & Gap\_BN & 7.79\% & \cellcolor[HTML]{D6D5D5} 2.29\% & \cellcolor[HTML]{D6D5D5} 1.24\% & \cellcolor[HTML]{D6D5D5} 1.16\% & \cellcolor[HTML]{D6D5D5} 1.19\% & \cellcolor[HTML]{D6D5D5} 1.5\% & \cellcolor[HTML]{D6D5D5} 1.71\% & \cellcolor[HTML]{D6D5D5} 0.41\% & \cellcolor[HTML]{D6D5D5} 0.02\% & \cellcolor[HTML]{E0FFFF} 0\% \\ \cline{2-13}
&	\multicolumn{1}{c|}{\multirow{6}{*}{0.4}} & M=LB & 0\% & 0\% & 0\% & 0\% & 0\% & 2\% & 3\% & 5\% & 5\% & 5\% \\ 
&	 & M$<$BN & \cellcolor[HTML]{D6D5D5} 100\% & 32\% & 6\% & 15\% & 1\% & 0\% & 0\% & 0\% & 0\% & 0\% \\ 
&	 & Gap\_M & \cellcolor[HTML]{D6D5D5} 1.82\% & 4.47\% & 6.21\% & 6.56\% & 3.93\% & 2.61\% & 2.55\% & 2.53\% & 2.52\% & 2.52\% \\ 
&	 & BN=LB & 0\% & 0\% & 0\% & 0\% & \cellcolor[HTML]{D6D5D5} 17\% & \cellcolor[HTML]{D6D5D5} 86\% & \cellcolor[HTML]{D6D5D5} 99\% & \cellcolor[HTML]{E0FFFF} 100\% & \cellcolor[HTML]{E0FFFF} 100\% & \cellcolor[HTML]{E0FFFF} 100\% \\ 
&	 & BN$<$M & 0\% & \cellcolor[HTML]{D6D5D5} 62\% & \cellcolor[HTML]{D6D5D5} 92\% & \cellcolor[HTML]{D6D5D5} 83\% & \cellcolor[HTML]{D6D5D5} 99\% & \cellcolor[HTML]{D6D5D5} 98\% & \cellcolor[HTML]{D6D5D5} 97\% & \cellcolor[HTML]{E0FFFF} 95\% & \cellcolor[HTML]{E0FFFF} 95\% & \cellcolor[HTML]{E0FFFF} 95\% \\ 
&	 & Gap\_BN & 8.57\% & \cellcolor[HTML]{D6D5D5} 3.51\% & \cellcolor[HTML]{D6D5D5} 2.67\% & \cellcolor[HTML]{D6D5D5} 3.17\% & \cellcolor[HTML]{D6D5D5} 1.19\% & \cellcolor[HTML]{D6D5D5} 0.07\% & \cellcolor[HTML]{D6D5D5} 0\% & \cellcolor[HTML]{E0FFFF} 0\% & \cellcolor[HTML]{E0FFFF} 0\% & \cellcolor[HTML]{E0FFFF} 0\% \\ \cline{2-13}
&	\multicolumn{1}{c|}{\multirow{6}{*}{0.5}} & M=LB & 0\% & 0\% & 0\% & 2\% & 12\% & 17\% & 17\% & 17\% & 17\% & 17\% \\ 
&	 & M$<$BN & \cellcolor[HTML]{D6D5D5} 100\% & 41\% & 46\% & 7\% & 3\% & 0\% & 0\% & 0\% & 0\% & 0\% \\ 
&	 & Gap\_M & \cellcolor[HTML]{D6D5D5} 1.95\% & 6.25\% & 4.93\% & 1.16\% & 0.92\% & 0.88\% & 0.88\% & 0.88\% & 0.88\% & 0.88\% \\ 
&	 & BN=LB & 0\% & 0\% & 0\% & \cellcolor[HTML]{D6D5D5} 39\% & \cellcolor[HTML]{D6D5D5} 94\% & \cellcolor[HTML]{E0FFFF} 100\% & \cellcolor[HTML]{E0FFFF} 100\% & \cellcolor[HTML]{E0FFFF} 100\% & \cellcolor[HTML]{E0FFFF} 100\% & \cellcolor[HTML]{E0FFFF} 100\% \\ 
&	 & BN$<$M & 0\% & \cellcolor[HTML]{D6D5D5} 55\% & \cellcolor[HTML]{D6D5D5} 48\% & \cellcolor[HTML]{D6D5D5} 84\% & \cellcolor[HTML]{D6D5D5} 88\% & \cellcolor[HTML]{E0FFFF} 83\% & \cellcolor[HTML]{E0FFFF} 83\% & \cellcolor[HTML]{E0FFFF} 83\% & \cellcolor[HTML]{E0FFFF} 83\% & \cellcolor[HTML]{E0FFFF} 83\% \\ 
&	 & Gap\_BN & 9.92\% & \cellcolor[HTML]{D6D5D5} 5.62\% & \cellcolor[HTML]{D6D5D5} 3.98\% & \cellcolor[HTML]{D6D5D5} 0.41\% & \cellcolor[HTML]{D6D5D5} 0.02\% & \cellcolor[HTML]{E0FFFF} 0\% & \cellcolor[HTML]{E0FFFF} 0\% & \cellcolor[HTML]{E0FFFF} 0\% & \cellcolor[HTML]{E0FFFF} 0\% & \cellcolor[HTML]{E0FFFF} 0\% \\ \cline{2-13}
&	\multicolumn{1}{c|}{\multirow{6}{*}{0.6}} & M=LB & 0\% & 0\% & 4\% & 24\% & 26\% & 26\% & 26\% & 26\% & 26\% & 26\% \\ 
&	 & M$<$BN & \cellcolor[HTML]{D6D5D5} 100\% & \cellcolor[HTML]{D6D5D5} 76\% & 36\% & 3\% & 0\% & 0\% & 0\% & 0\% & 0\% & 0\% \\ 
&	 & Gap\_M & \cellcolor[HTML]{D6D5D5} 2.37\% & \cellcolor[HTML]{D6D5D5} 5.02\% & 0.79\% & 0.44\% & 0.43\% & 0.43\% & 0.43\% & 0.43\% & 0.43\% & 0.43\% \\ 
&	 & BN=LB & 0\% & 0\% & \cellcolor[HTML]{D6D5D5} 6\% & \cellcolor[HTML]{D6D5D5} 97\% & \cellcolor[HTML]{E0FFFF} 100\% & \cellcolor[HTML]{E0FFFF} 100\% & \cellcolor[HTML]{E0FFFF} 100\% & \cellcolor[HTML]{E0FFFF} 100\% & \cellcolor[HTML]{E0FFFF} 100\% & \cellcolor[HTML]{E0FFFF} 100\% \\ 
&	 & BN$<$M & 0\% & 21\% & \cellcolor[HTML]{D6D5D5} 46\% & \cellcolor[HTML]{D6D5D5} 76\% & \cellcolor[HTML]{E0FFFF} 74\% & \cellcolor[HTML]{E0FFFF} 74\% & \cellcolor[HTML]{E0FFFF} 74\% & \cellcolor[HTML]{E0FFFF} 74\% & \cellcolor[HTML]{E0FFFF} 74\% & \cellcolor[HTML]{E0FFFF} 74\% \\ 
&	 & Gap\_BN & 11.83\% & 7.04\% & \cellcolor[HTML]{D6D5D5} 0.65\% & \cellcolor[HTML]{D6D5D5} 0.01\% & \cellcolor[HTML]{E0FFFF} 0\% & \cellcolor[HTML]{E0FFFF} 0\% & \cellcolor[HTML]{E0FFFF} 0\% & \cellcolor[HTML]{E0FFFF} 0\% & \cellcolor[HTML]{E0FFFF} 0\% & \cellcolor[HTML]{E0FFFF} 0\% \\ \cline{2-13}
&	\multicolumn{1}{c|}{\multirow{6}{*}{0.7}} & M=LB & 0\% & \cellcolor[HTML]{D6D5D5} 1\% & 32\% & 62\% & 62\% & 62\% & 62\% & 62\% & 62\% & 62\% \\ 
&	 & M$<$BN & \cellcolor[HTML]{D6D5D5} 100\% & \cellcolor[HTML]{D6D5D5} 100\% & 28\% & 0\% & 0\% & 0\% & 0\% & 0\% & 0\% & 0\% \\ 
&	 & Gap\_M & \cellcolor[HTML]{D6D5D5} 5.53\% & \cellcolor[HTML]{D6D5D5} 0.85\% & 0.22\% & 0.14\% & 0.14\% & 0.14\% & 0.14\% & 0.14\% & 0.14\% & 0.14\% \\ 
&	 & BN=LB & 0\% & 0\% & \cellcolor[HTML]{D6D5D5} 47\% & \cellcolor[HTML]{E0FFFF} 100\% & \cellcolor[HTML]{E0FFFF} 100\% & \cellcolor[HTML]{E0FFFF} 100\% & \cellcolor[HTML]{E0FFFF} 100\% & \cellcolor[HTML]{E0FFFF} 100\% & \cellcolor[HTML]{E0FFFF} 100\% & \cellcolor[HTML]{E0FFFF} 100\% \\
&	 & BN$<$M & 0\% & 0\% & \cellcolor[HTML]{D6D5D5} 46\% & \cellcolor[HTML]{E0FFFF} 38\% & \cellcolor[HTML]{E0FFFF} 38\% & \cellcolor[HTML]{E0FFFF} 38\% & \cellcolor[HTML]{E0FFFF} 38\% & \cellcolor[HTML]{E0FFFF} 38\% & \cellcolor[HTML]{E0FFFF} 38\% & \cellcolor[HTML]{E0FFFF} 38\% \\ 
&	 & Gap\_BN & 14.64\% & 2.81\% & \cellcolor[HTML]{D6D5D5} 0.15\% & \cellcolor[HTML]{E0FFFF} 0\% & \cellcolor[HTML]{E0FFFF} 0\% & \cellcolor[HTML]{E0FFFF} 0\% & \cellcolor[HTML]{E0FFFF} 0\% & \cellcolor[HTML]{E0FFFF} 0\% & \cellcolor[HTML]{E0FFFF} 0\% & \cellcolor[HTML]{E0FFFF} 0\% \\ \cline{2-13}
&	\multicolumn{1}{c|}{\multirow{6}{*}{0.8}} & M=LB & 0\% & \cellcolor[HTML]{D6D5D5} 3\% & 71\% & 98\% & 98\% & 98\% & 98\% & 98\% & 98\% & 98\% \\ 
&	 & M$<$BN & \cellcolor[HTML]{D6D5D5} 100\% & \cellcolor[HTML]{D6D5D5} 92\% & 7\% & 0\% & 0\% & 0\% & 0\% & 0\% & 0\% & 0\% \\ 
&	 & Gap\_M & \cellcolor[HTML]{D6D5D5} 3.78\% & \cellcolor[HTML]{D6D5D5} 0.42\% & 0.06\% & 0\% & 0\% & 0\% & 0\% & 0\% & 0\% & 0\% \\ 
&	 & BN=LB & 0\% & 0\% & \cellcolor[HTML]{D6D5D5} 87\% & \cellcolor[HTML]{E0FFFF} 100\% & \cellcolor[HTML]{E0FFFF} 100\% & \cellcolor[HTML]{E0FFFF} 100\% & \cellcolor[HTML]{E0FFFF} 100\% & \cellcolor[HTML]{E0FFFF} 100\% & \cellcolor[HTML]{E0FFFF} 100\% & \cellcolor[HTML]{E0FFFF} 100\% \\ 
&	 & BN$<$M & 0\% & 4\% & \cellcolor[HTML]{D6D5D5} 23\% & \cellcolor[HTML]{E0FFFF} 2\% & \cellcolor[HTML]{E0FFFF} 2\% & \cellcolor[HTML]{E0FFFF} 2\% & \cellcolor[HTML]{E0FFFF} 2\% & \cellcolor[HTML]{E0FFFF} 2\% & \cellcolor[HTML]{E0FFFF} 2\% & \cellcolor[HTML]{E0FFFF} 2\% \\ 
&	 & Gap\_BN & 9\% & 1.14\% & \cellcolor[HTML]{D6D5D5} 0.02\% & \cellcolor[HTML]{E0FFFF} 0\% & \cellcolor[HTML]{E0FFFF} 0\% & \cellcolor[HTML]{E0FFFF} 0\% & \cellcolor[HTML]{E0FFFF} 0\% & \cellcolor[HTML]{E0FFFF} 0\% & \cellcolor[HTML]{E0FFFF} 0\% & \cellcolor[HTML]{E0FFFF} 0\% \\ \cline{2-13}
&	\multicolumn{1}{c|}{\multirow{6}{*}{0.9}} & M=LB & 0\% & \cellcolor[HTML]{D6D5D5} 14\% & \cellcolor[HTML]{D6D5D5} 85\% & \cellcolor[HTML]{E0FFFF} 100\% & \cellcolor[HTML]{E0FFFF} 100\% & \cellcolor[HTML]{E0FFFF} 100\% & \cellcolor[HTML]{E0FFFF} 100\% & \cellcolor[HTML]{E0FFFF} 100\% & \cellcolor[HTML]{E0FFFF} 100\% & \cellcolor[HTML]{E0FFFF} 100\% \\
&	 & M$<$BN & \cellcolor[HTML]{D6D5D5} 100\% & \cellcolor[HTML]{D6D5D5} 87\% & \cellcolor[HTML]{D6D5D5} 12\% & \cellcolor[HTML]{E0FFFF} 0\% & \cellcolor[HTML]{E0FFFF} 0\% & \cellcolor[HTML]{E0FFFF} 0\% & \cellcolor[HTML]{E0FFFF} 0\% & \cellcolor[HTML]{E0FFFF} 0\% & \cellcolor[HTML]{E0FFFF} 0\% & \cellcolor[HTML]{E0FFFF} 0\% \\ 
&	 & Gap\_M & \cellcolor[HTML]{D6D5D5} 2.61\% & \cellcolor[HTML]{D6D5D5} 0.25\% & \cellcolor[HTML]{D6D5D5} 0.02\% & \cellcolor[HTML]{E0FFFF} 0\% & \cellcolor[HTML]{E0FFFF} 0\% & \cellcolor[HTML]{E0FFFF} 0\% & \cellcolor[HTML]{E0FFFF} 0\% & \cellcolor[HTML]{E0FFFF} 0\% & \cellcolor[HTML]{E0FFFF} 0\% & \cellcolor[HTML]{E0FFFF} 0\% \\ 
&	 & BN=LB & 0\% & 1\% & 82\% & \cellcolor[HTML]{E0FFFF} 100\% & \cellcolor[HTML]{E0FFFF} 100\% & \cellcolor[HTML]{E0FFFF} 100\% & \cellcolor[HTML]{E0FFFF} 100\% & \cellcolor[HTML]{E0FFFF} 100\% & \cellcolor[HTML]{E0FFFF} 100\% & \cellcolor[HTML]{E0FFFF} 100\% \\ 
&	 & BN$<$M & 0\% & 5\% & 8\% & \cellcolor[HTML]{E0FFFF} 0\% & \cellcolor[HTML]{E0FFFF} 0\% & \cellcolor[HTML]{E0FFFF} 0\% & \cellcolor[HTML]{E0FFFF} 0\% & \cellcolor[HTML]{E0FFFF} 0\% & \cellcolor[HTML]{E0FFFF} 0\% & \cellcolor[HTML]{E0FFFF} 0\% \\ 
&	 & Gap\_BN & 5.34\% & 0.73\% & 0.03\% & \cellcolor[HTML]{E0FFFF} 0\% & \cellcolor[HTML]{E0FFFF} 0\% & \cellcolor[HTML]{E0FFFF} 0\% & \cellcolor[HTML]{E0FFFF} 0\% & \cellcolor[HTML]{E0FFFF} 0\% & \cellcolor[HTML]{E0FFFF} 0\% & \cellcolor[HTML]{E0FFFF} 0\% \\ \hline
\multicolumn{1}{c}{\multirow{1}{*}{}} & & tmin\_BN & 0.3095 & 0.1534 & 0.0907 & 0.0292 & 0.0294 & 0.0284 & 0.0272 & 0.0275 & 0.0255 & 0.0255\\
\multicolumn{1}{c}{\multirow{1}{*}{}} & & tmax\_BN & 5.0831 & 2.2143 & 1.1515 & 0.7193 & 0.4101 & 0.2772 & 0.2165 & 0.1726 & 
0.1419 & 0.1235\\
\multicolumn{1}{c}{\multirow{1}{*}{}} & & tavg\_BN & 1.8845 & 0.7478 &	0.367 & 0.2151 & 0.1319 & 0.095 & 0.0762 & 0.0638 & 0.0549 & 0.049\\ \cline{3-13}

\end{tabular}
}

\caption{Computational results obtained by algorithms $M$  and algorithm $BN$ on $TI(1000,B,\Delta)$}
\label{BNvsM1000TI}

\end{center}

\end{table}

We also measured the average time (in seconds) required by $BN$ to solve one instance out of the 100 in each cell. At the bottom of each column of Table \ref{BNvsM1000TI} we report the minimum, maximum, and average time of these values (tmin\_BN, tmax\_BN, tavg\_BN, respectively). We remark that the maximum is always reached for $\Delta = 0.1$. Basically these times decrease for increasing $B$: we think that this is due to the reduced number of operations in {\sc Phase II}. The average times required by algorithms $M$ to solve one instance out of the 100 in each cell are not displayed because they are always smaller than 0.1 seconds for all $B$ and $\Delta$. 

In Table \ref{BNvsMTot120to1000TI}, the results for $n \in \{120,$ $250,$ $500,$ $1000\}$ are summarized (in column \linebreak $B \in \{120,150\}$ they are averaged over the corresponding 2000 instances, in column \linebreak $B = 180$ over 1000, in column $B \in \{210,\dots,390\}$ over 7000, and in column \linebreak $B \in \{120,\dots,390\}$ over 10000). Results essentially reflect what happens for $n = 1000$. In particular, algorithm $BN$ is definitely better than algorithms $M$ for $B \ge 210$, algorithms $M$ have better performances for $B \le 150$, while for $B = 180$ algorithm $BN$ improves with the growth of $n$. We also remark that for $B \in \{120,150\}$ the values M=LB and BN=LB are very small: we suspect that this is due to the poor quality of $LB_{BPPC}$. In fact, the experiments conducted in \cite{BN2017IGexact} with $n = 250$ show that $LB$ is strictly smaller than the value of the optimum solution in 98 out of 100 instances with $B = 120$, 23 out of 100 for $B = 150$, 3 out of 100 for $B = 180$, zero in all the other cases.

\begin{table}

\begin{center}

{\scriptsize
\begin{tabular}{| c | c | r | c | c | c || c | c | }
\cline{4-8}
\multicolumn{1}{c}{} & \multicolumn{1}{c}{} & \multicolumn{1}{c|}{}                     &                  \multicolumn{4}{c|}{$B$} \\                                                                                                                                                                                                                                                                                                                                                                                                              

\cline{4-8}
\multicolumn{1}{c}{} & \multicolumn{1}{c}{} & \multicolumn{1}{c|}{}  & $\{120.150\}$ & $180$ & $\{210.\dots.390\}$ & $\{120.\dots.390\}$ \\ \hline
\multicolumn{1}{|c|}{\multirow{24}{*}{$n$}} & \multicolumn{1}{c|}{\multirow{6}{*}{120}} & M=LB  & \cellcolor[HTML]{D6D5D5} 11.95\% & \cellcolor[HTML]{D6D5D5} 44.4\% & 65.84\% & 52.92\% \\ 
& & M$<$BN  & \cellcolor[HTML]{D6D5D5} 83.6\% & \cellcolor[HTML]{D6D5D5} 18.1\% & 1.2\% & \cellcolor[HTML]{D6D5D5}19.37\% \\ 
& & Gap\_M  & \cellcolor[HTML]{D6D5D5} 4.91\% & \cellcolor[HTML]{D6D5D5} 2.23\% & 1.52\%  & \cellcolor[HTML]{D6D5D5} 2.27\% \\
& & BN=LB  & 4.4\% & 42.3\% & \cellcolor[HTML]{D6D5D5} 83.04\% & \cellcolor[HTML]{D6D5D5} 63.24\% \\ 
& & BN$<$M  & 1.25\% & 13.8\% & \cellcolor[HTML]{D6D5D5} 22.73\% & 17.54\% \\
& & Gap\_BN  & 9.26\% & 2.35\% & \cellcolor[HTML]{D6D5D5} 0.68\% & 2.56\% \\ \cline{2-8}

& \multicolumn{1}{c|}{\multirow{6}{*}{250}} & M=LB  & \cellcolor[HTML]{D6D5D5} 6.15\% & \cellcolor[HTML]{D6D5D5} 33.8\% & 47.63\% & 37.95\% \\ 
& & M$<$BN  & \cellcolor[HTML]{D6D5D5} 89.45\% & 17.7\% & 0.54\% & 20.04\% \\ 
& & Gap\_M  & \cellcolor[HTML]{D6D5D5} 4.07\% & 2.07\% & 1.72\% & 2.23\% \\
& & BN=LB  & 1.3\% & 32.3\% & \cellcolor[HTML]{D6D5D5} 79.14\% & \cellcolor[HTML]{D6D5D5} 58.89\% \\ 
& & BN$<$M  & 2.95\% & \cellcolor[HTML]{D6D5D5} 26.2\% & \cellcolor[HTML]{D6D5D5} 47.41\% & \cellcolor[HTML]{D6D5D5} 36.4\% \\
& & Gap\_BN  & 8.18\% & \cellcolor[HTML]{D6D5D5} 1.9\% & \cellcolor[HTML]{D6D5D5} 0.41\% & \cellcolor[HTML]{D6D5D5} 2.16\% \\ \cline{2-8} 

& \multicolumn{1}{c|}{\multirow{6}{*}{500}} & M=LB  & \cellcolor[HTML]{D6D5D5} 3.75\% & 27.2\% & 37.64\% & 29.82\% \\ 
& & M$<$BN  & \cellcolor[HTML]{D6D5D5} 85.4\% & 15.1\% & 0.27\% & 18.78\% \\ 
& & Gap\_M  & \cellcolor[HTML]{D6D5D5} 2.99\% & 1.92\% & 1.87\% & 2.1\% \\
& & BN=LB  & 0.75\% & \cellcolor[HTML]{D6D5D5} 27.8\% & \cellcolor[HTML]{D6D5D5} 75.87\% & \cellcolor[HTML]{D6D5D5} 56.04\% \\ 
& & BN$<$M  & 6.8\% & \cellcolor[HTML]{D6D5D5} 46.6\% & \cellcolor[HTML]{D6D5D5} 61.47\% & \cellcolor[HTML]{D6D5D5} 49.05\% \\
& & Gap\_BN  & 6.85\% & \cellcolor[HTML]{D6D5D5} 1.45\% & \cellcolor[HTML]{D6D5D5} 0.35\% & \cellcolor[HTML]{D6D5D5} 1.76\% \\ \cline{2-8}

&  \multicolumn{1}{c|}{\multirow{6}{*}{1000}} & M=LB  & \cellcolor[HTML]{D6D5D5} 0.9\% & 9.6\% & 30.27\% & 23.29\% \\ 
& & M$<$BN  & \cellcolor[HTML]{D6D5D5} 80.1\% & 6.8\% & 0.41\% & 17.67\% \\ 
& & Gap\_M  & \cellcolor[HTML]{D6D5D5} 2.44\% & 0.98\% & 2.05\% & 2.12\% \\
& & BN=LB  & 0.05\% & \cellcolor[HTML]{D6D5D5} 11.1\% & \cellcolor[HTML]{D6D5D5} 73.89\% & \cellcolor[HTML]{D6D5D5} 53.95\% \\ 
& & BN$<$M  & 14.1\% & \cellcolor[HTML]{D6D5D5} 32.4\% & \cellcolor[HTML]{D6D5D5} 69.26\% & \cellcolor[HTML]{D6D5D5} 57.78\% \\
& & Gap\_BN  & 5.93\% & \cellcolor[HTML]{D6D5D5} 0.52\% & \cellcolor[HTML]{D6D5D5} 0.24\% & \cellcolor[HTML]{D6D5D5} 1.46\% \\ \hline
	
\end{tabular}
}

\caption{General results obtained by algorithms $M$  and algorithm $BN$ on $TI(n,B,\Delta)$}
\label{BNvsMTot120to1000TI}

\end{center}

\end{table}
\section{Computational results on literature instances}
\label{sec:exactresultsTM}

In this section we  discuss the results obtained by solving the instances by \cite{MIMT2010}  (see http://www.or.deis.unibo.it) and other instances with threshold conflict graphs by running algorithms $M$ and $BN$ on them.

A graph is a threshold graph if there exist a real number $d$ (the threshold) and a weight $p_x$ for every vertex $x$ such that $(i,j)$ is an edge iff $(p_i + p_j)/2 \le d$  (\cite{CH1973}). A threshold graph has many peculiar properties as it is at the same time an interval graph, a co-interval graph, a cograph, a split graph, and a permutation graph (\cite{G1980}). In addition, its complement, where $(i,j)$ is an edge iff $(p_i + p_j)/2 > d$, is a threshold graph too. $VC$ is solvable in linear time on threshold graphs, too, nevertheless $BPPC$ with a threshold conflict graph remains $NP$-hard.

\cite{GLS2004} describe the following generator, which we shall refer to as {\em {\sc t}-generator}:
``{\em A value $p_i$ was first assigned to each vertex $i \in V$ according to a continuous uniform distribution on $[0,1]$. Each edge $(i,j)$ of $G$ was created whenever $(p_i + p_j)/2 \le d$, where $d$ is the expected density of $G$.}''  This generator clearly produces threshold graphs. In addition, the expected edge density $\delta$ is not $d$ as claimed. Actually, $\delta$ is a function of $d$ (\cite{BN2017TG}), precisely: $\delta = f(d) = \frac{2(nd)^2 - nd}{n(n-1)}$ for $d \le 0.5$ and $\delta = f(d)= \frac{n(n-1) - 2n^2(1-d)^2 - n(1-d)}{n(n-1)}$ for $d \ge 0.5$; here we recall some useful pairs $(d,\delta)$: (0,0), (0.1,0.02), (0.2,0.08), (0.3,0.18), (0.4,0.32), (0.5,0.5), (0.6,0.68), (0.7,0.82), (0.8,0.92), (0.9,0.98), (1,1).

This generator has been improperly used to generate arbitrary graphs, and, in particular, \cite{MIMT2010} made publicly available (http://www.or.deis.unibo.it) instances generated in this way and used by many authors (see list in Section \ref{sec:literature}). \cite{SV2013} observed that the {\sc t}-generator generates interval conflict graphs (and, in fact, their Dynamic Programming phase is designed for interval graphs), but actually these graphs are special interval graphs.

By $TM(n,B,f(d))$ we denote a set of ten instances with $n$ items, bound $B$, and threshold conflict graph with density $f(d)$. In particular $n \in \{120,250,500,1000\}$, \linebreak and $d \in \{0,0.1,\dots,0.9\}$. The weights and the conflict graphs of all the $TM(n,B,f(d))$ are exactly those in the classes 1,2,3,4 by \cite{MIMT2010}. As for $B$, we considered $B \in \{120,150,\dots,390\} \cup \{400\}$, even if in the cited paper only $B=150$ is considered. In particular, \cite{MIMT2010} select the first 10 instances of the 20 originally proposed by \cite{F96} for the Bin Packing (without conflicts), and add 10 random threshold conflict graphs generated by means of the {\sc t}-generator, varying $d$ from 0 to 0.9. The Bin Packing instances proposed by \cite{F96} have weights uniformly distributed in $[20,100]$ and $B=150$ because, as the author says, this setup was the most difficult for the Bin Packing lower bound algorithms by \cite{MT1990}. Nevertheless, \cite{G1998} easily solves the last five open instances. The instances $TM(n,150,f(d))$ correspond to the so-called ``$u$ instances'' by  \cite{SV2013}.

In order to verify how much the item weights affect the quality of the solution and/or the computing time, we also decided to construct the $TS$ instances: by $TS(n,B,f(d))$ we will denote a set of ten instances with $n$ items, bound $B$, and threshold conflict graph with density $f(d)$. The conflict graphs of a $TS(n,\cdot,f(d))$ are those of $TM(n,\cdot,f(d))$, and the weights  are uniformly distributed in $[500,2500]$. We choose $B \in \{3000, 3750, \dots, 9750\} \cup \{10000\}$. We remark that the item weights of $TS(n,B,f(d))$ are generated as the ``instances with a larger number of items per bin'' by  \cite{SV2013} (the  so-called ``$d$ instances''), where, however, only $B = 10000$ is considered.

Let $\overline{w}$ be the average weight of an item, then the average number of items per bin is $B / \overline{w}$. It is worth observing that the same number of items per bin is obtained in the $TM(n,B,f(d))$ and in the $TS(n,25 \times B,f(d))$ instances.

Observe also that in the $TS$ instances the number of different weights is $2500-500+1=2001$. Hence, every weight is expected to appear in $n/2001$ copies. For a same $n$, the (classical) Bin Packing underlying a $TM$ instance recalls a Cutting Stock problem (see Section \ref{sec:resultsIG}) in a stronger way than the one underlying a $TS$ instance.

\vspace{0.5cm}
\noindent
In Table \ref{BNvsMTot500to1000TM} we compare the results obtained by the heuristic algorithms $M$ and $BN$ over $TM(500,B,f(d))$ and $TM(1000,B,f(d))$ (in column $B \in \{120,150\}$ the results are averaged over 2000 instances, in column $B = 180$ over 1000, in column $B \in \{210,\dots,390\} \cup \{400\}$ over 8000, and in column $B \in \{120,\dots,390\} \cup \{400\}$ over 11000). 

The results of the heuristic algorithms $M$ and $BN$ over $TS(500,B,f(d))$ and \linebreak $TS(1000,B,f(d))$ can be found in Table \ref{BNvsMTot500to1000TS} (in column $B \in \{3000,3750\}$ the results are averaged over 2000 instances, in column $B = 4500$ over 1000, in column $B \in \{5250,\dots,9750\} \cup \{10000\}$ over 8000, and in column $B \in \{3000,\dots,9750\} \cup \{10000\}$ over 11000).  

We remind that the values in an arbitrary column in Table \ref{BNvsMTot500to1000TM} (Table \ref{BNvsMTot500to1000TS}, respectively) are the average of the values obtained on $TM(n,B,f(d))$ ($TS(n,B,f(d))$, respectively) for the corresponding $B$'s and $d=0,0.1,\dots,0.9$. 

The results in Tables \ref{BNvsMTot500to1000TM} and \ref{BNvsMTot500to1000TS} show that the behaviour of each heuristic algorithm w.r.t.~the quality of the solution on the instances  $TM(n,B,f(d))$ and $TS(n,B,f(d))$ is essentially the same. We can note that $BN$ outperforms $M$ for $B \ge 180$  on all the $TM$ instances and for the corresponding $B \ge 4500$ on all the $TS$ instances, i.e.~when the average number of items per bin is greater than or equal to 3.5.

\begin{table}

\begin{center}

{\scriptsize
\begin{tabular}{ | c | c | r | c | c | c || c | }

\cline{4-7}
\multicolumn{1}{c}{} & \multicolumn{1}{c}{} & \multicolumn{1}{c|}{}                     &                  \multicolumn{4}{c|}{$B$} \\                                                                                                                                                                                                                                                                                                                                                                                                              

\cline{4-7}
\multicolumn{1}{c}{} & \multicolumn{1}{c}{} & \multicolumn{1}{c|}{}  & $\{120,150\}$ & $180$ & $\{210,\dots,390\}$ & $\{120,\dots,390\} \cup \{400\}$ \\ \hline
\multicolumn{1}{|c|}{\multirow{12}{*}{$n$}} & \multicolumn{1}{c|}{\multirow{6}{*}{500}} & tM=LB  & \cellcolor[HTML]{D6D5D5} 3\% & 28\% & 78.57\% & 58.4\% \\ 
& & tM$<$tBN  & \cellcolor[HTML]{D6D5D5} 85\% & 9\% & 0.71\% & 18.4\% \\ 
& & Gap\_tM  & \cellcolor[HTML]{D6D5D5} 3.18\% & 1.46\% & 0.76\% & \cellcolor[HTML]{D6D5D5} 1.32\% \\
& & tBN=LB  & 0\% & \cellcolor[HTML]{D6D5D5} 31\% & \cellcolor[HTML]{D6D5D5} 93\% & \cellcolor[HTML]{D6D5D5} 68.2\% \\ 
& & tBN$<$tM  & 7\% & \cellcolor[HTML]{D6D5D5} 47\% & \cellcolor[HTML]{D6D5D5} 28.57\% & \cellcolor[HTML]{D6D5D5} 26.1\% \\
& & Gap\_tBN  & 5.97\% & \cellcolor[HTML]{D6D5D5} 0.61\% & \cellcolor[HTML]{D6D5D5} 0.16\% & 1.37\% \\ \cline{2-7}

&  \multicolumn{1}{c|}{\multirow{6}{*}{1000}} & tM=LB  & \cellcolor[HTML]{D6D5D5} 0.5\% & 18\% & 77.43\% & 56.1\% \\ 
& & tM$<$tBN  & \cellcolor[HTML]{D6D5D5} 82.5\% & 4\% & 1\% & 17.6\% \\ 
& & Gap\_tM  & \cellcolor[HTML]{D6D5D5} 2.29\% & 1.36\% & 0.69\% & 1.08\% \\
& & tBN=LB  & 0\% & \cellcolor[HTML]{D6D5D5} 29\% & \cellcolor[HTML]{D6D5D5} 97.29\% & \cellcolor[HTML]{D6D5D5} 71\% \\ 
& & tBN$<$tM  & 11.5\% & \cellcolor[HTML]{D6D5D5} 59\% & \cellcolor[HTML]{D6D5D5} 28.71\% & \cellcolor[HTML]{D6D5D5} 28.3\% \\
& & Gap\_tBN  & 4.96\% & \cellcolor[HTML]{D6D5D5} 0.39\% & \cellcolor[HTML]{D6D5D5} 0.06\% & \cellcolor[HTML]{D6D5D5} 1.07\% \\ \hline
	
\end{tabular}
}

\caption{General results obtained by algorithms $M$  and algorithm $BN$ on $TM(n,B,f(d))$}
\label{BNvsMTot500to1000TM}

\end{center}

\end{table}

\begin{table}

\begin{center}

{\scriptsize
\begin{tabular}{| c | c | r | c | c | c || c | }

\cline{4-7}
\multicolumn{1}{c}{} & \multicolumn{1}{c}{} & \multicolumn{1}{c|}{} &                  \multicolumn{4}{c|}{$B$} \\                                                                                                                                                                                                                                                                                                                                                                                                              

\cline{4-7}
\multicolumn{1}{c}{} & \multicolumn{1}{c}{} & \multicolumn{1}{c|}{}  & $\{3000,3750\}$ & $4500$ & $\{5250,\dots,9750\}$ & $\{3000,\dots,9750\} \cup \{10000\}$ \\ \hline
\multicolumn{1}{|c|}{\multirow{12}{*}{$n$}} & \multicolumn{1}{c|}{\multirow{6}{*}{500}} & sM=LB  & \cellcolor[HTML]{D6D5D5} 3\% & 31\% & 78.86\% & 58.9\% \\ 
& & sM$<$sBN  & \cellcolor[HTML]{D6D5D5} 83.5\% & 9\% & 0.86\% & 18.2\% \\ 
& & Gap\_sM  & \cellcolor[HTML]{D6D5D5} 3.01\% & 1.52\% & 0.78\% & \cellcolor[HTML]{D6D5D5} 1.3\% \\
& & sBN=LB  & 0.5\% & \cellcolor[HTML]{D6D5D5} 34\% & \cellcolor[HTML]{D6D5D5} 95.86\% & \cellcolor[HTML]{D6D5D5} 70.6\% \\ 
& & sBN$<$sM  & 9\% & \cellcolor[HTML]{D6D5D5} 52\% & \cellcolor[HTML]{D6D5D5} 28.86\% & \cellcolor[HTML]{D6D5D5} 27.2\% \\
& & Gap\_sBN  & 5.78\% & \cellcolor[HTML]{D6D5D5} 0.61\% & \cellcolor[HTML]{D6D5D5} 0.15\% & 1.32\% \\ \cline{2-7}

&  \multicolumn{1}{c|}{\multirow{6}{*}{1000}} & sM=LB  & \cellcolor[HTML]{D6D5D5} 1\% & 24\% & 77.14\% & 56.6\% \\ 
& & sM$<$sBN  & \cellcolor[HTML]{D6D5D5} 87\% & 5\% & 1\% & 18.6\% \\ 
& & Gap\_sM  & \cellcolor[HTML]{D6D5D5} 2.31\% & 1.46\% & 0.73\% & \cellcolor[HTML]{D6D5D5} 1.12\% \\
& & sBN=LB  & 0\% & \cellcolor[HTML]{D6D5D5} 28\% & \cellcolor[HTML]{D6D5D5} 90.43\% & \cellcolor[HTML]{D6D5D5} 66.1\% \\ 
& & sBN$<$sM  & 7.5\% & \cellcolor[HTML]{D6D5D5} 59\% & \cellcolor[HTML]{D6D5D5} 28.29\% & \cellcolor[HTML]{D6D5D5} 27.2\% \\
& & Gap\_sBN  & 5.13\% & \cellcolor[HTML]{D6D5D5} 0.46\% & \cellcolor[HTML]{D6D5D5} 0.09\% & 1.14\% \\ \hline
	
\end{tabular}
}

\caption{General results obtained by algorithms $M$  and algorithm $BN$ on $TS(n,B,f(d))$}
\label{BNvsMTot500to1000TS}

\end{center}

\end{table}

\section{Conclusions}
\label{sec:conclusions}

\noindent
In this paper we dealt with the Bin Packing Problem with Conflicts ($BPPC$) on instances with interval conflict graphs.

We proposed a new heuristic algorithm for the problem. We conducted experiments by varying the number $n$ of items, the edge density of the conflict graph, the value $B$, and the values of the weights, hence the average number of items per bin. We remark that the experiments in the papers \cite{BW2005, BCM2013, BP2016, ELGN2011, GI2016, KCHT2012, KCT2010, KCT2012, MG2009, MR2011, MIMT2010, YLW2014}  consider weights uniformly distributed in $[20;100]$, $B=150$, and threshold conflict graphs, only.

We compare the results of our algorithm to the results obtained by running a parameterized adaptation of three classical heuristic algorithms for $BP$. The results show that our algorithm outperforms them definitively when the average number of items per bin is greater than or equal to 3.5, both on instances with interval and threshold conflict graphs.   

To our knowledge, no random interval graph generator exists which outputs a graph with desired edge density. For this reason, we defined a new one with this properties, and used it to generate thousands of instances on which, grouped by edge density, we tested our algorithm.   


\end{document}